\documentclass[a4paper,english]{smfart}

\usepackage{smfthm}
\usepackage{xypic}

\usepackage{pstricks}

\newcommand{\SO}{\mbox{\rm SO}}
\newcommand{\Spin}{\mbox{\rm Spin}}
\newcommand{\PSOM}{P_{\rm\scriptsize \SO}M}
\newcommand{\Aut}{\mbox{\rm Aut}}
\newcommand{\Cl}{\mbox{\rm Cl}}                
\newcommand{\Res}{\mbox{\rm Res}}                
\newcommand{\sign}{\mbox{\rm sign}}
\newcommand{\sys}{\mbox{\rm sys}}
\newcommand{\spinsys}{\mbox{\rm spin-sys}}
\newcommand{\length}{\mbox{\rm length}}
\newcommand{\area}{\mbox{\rm area}}
\newcommand{\vol}{\mbox{\rm vol}}
\newcommand{\<}{\left\langle}
\renewcommand{\>}{\right\rangle}
\newcommand{\E}{\mathcal E}
\newcommand{\CC}{\mathbb C}
\newcommand{\NN}{\mathbb N}
\newcommand{\RR}{\mathbb R}
\newcommand{\ZZ}{\mathbb Z}
\newcommand{\PP}{\mathbb P}

\author{Christian B\"ar}

\address{Universit\"at Hamburg\\ FB Mathematik\\ Bundesstr.~55\\ 
         D-20146 Hamburg}
\email{baer{@}math.uni-hamburg.de}
\urladdr{http://www.math.uni-hamburg.de/home/baer/}

\title{Dependence on the spin structure of the Dirac spectrum}

\begin{document}

\begin{abstract}
The theme is the influence of the spin structure on the Dirac
spectrum of a spin manifold.
We survey examples and results related to this question.
\end{abstract}

\alttitle{D\'ependance du Spectre de l'Op\'erateur de Dirac de la 
Structure Spinorielle}

\begin{altabstract}
Sur une vari\'et\'e spinorielle, nous \'etudions la d\'ependance du spectre de 
l'op\'erateur de Dirac par rapport \`a la structure spinorielle. 
Nous donnons un r\'esum\'e des exemples et des r\'esultats li\'es 
\`a cette question.
\end{altabstract}


\maketitle

\tableofcontents

\section{Introduction}

The relation between the geometry of a Riemannian manifold and the
spectrum of its Laplace operator acting on functions (or more generally,
on differential forms), has attracted a lot of attention.
This is the question how shape and sound of a space are related.
A beautiful introduction into this topic can be found in \cite{chavel84a}.
When one passes from this ``bosonic'' theory to ``fermions'', i.e.\ when
turning to spinors and the Dirac operator, a new object enters the stage,
the {\em spin structure}.
This is a global topological object needed to define spinors.
The question arises how this piece of structure, in addition to the usual 
geometry of the manifold, influences the spectrum of the Dirac operator.

It has been known for a long time that even on the simplest examples such as 
the 1-sphere the Dirac spectrum does depend on the spin structure.
We will discuss the 1-sphere, flat tori, 3-dimensional Bieberbach manifolds,
and spherical space forms in some detail.
For these manifolds the spectrum can be computed explicitly.
For some of these examples an important invariant computed out of the
spectrum, the $\eta$-invariant, also depends on the spin structure.
On the other hand, under a certain assumption, the difference between
the $\eta$-invariants for two spin structures on the same manifold must be 
an integer. 
Hence the two $\eta$-invariants are not totally unrelated. 

We also look at circle bundles and the behavior of the Dirac spectrum
under collapse.
This means that one shrinks the fibers to points.
The spin structure determines the qualitative spectral behavior.
If the spin structure is projectable, then some eigenvalues tend to 
$\pm\infty$ while the others essentially converge to the eigenvalues 
of the basis manifold.
If the spin structure if nonprojectable, then all eigenvalues diverge.

In most examples it is totally hopeless to try to explicitly compute the
Dirac (or other) spectra.
Still, eigenvalue estimates are very often possible.
So far, these estimates have not taken into account the spin structure
despite its influence on the spectrum.
The reason for this lies in the essentially local methods such as the
Bochner technique.
In order to get better estimates taking the spin structure into account
one first has to find new, truly spin geometric invariants.
We discuss some of the first steps in this direction.
Here the {\em spinning systole} is the relevant spin geometric input.

Finally we look at noncompact examples in order to check if the continuous
spectrum is affected by a change of spin structure.
It turns out that this is the case.
There are hyperbolic manifolds having two spin structures such that for
the first one the Dirac spectrum is discrete while it is all of $\RR$ for
the other one.
The influence of the spin structure could hardly be any more dramatic.

{\em Acknowledgements.}
It is a pleasure to thank B.~Ammann, M.~Dahl, and F.~Pf\"affle for helpful
discussion.

\section{Generalities}
Let us start by collecting some terminology and basic facts.
A more thorough introduction to the concepts of spin geometry
can e.g.\ be found in \cite{gilkey84a,berline-getzler-vergne91a,
lawson-michelsohn89a}.
Let $M$ denote an $n$-dimensional oriented Riemannian manifold with a spin 
structure $P$.
This is a $\Spin(n)$-principal bundle which doubly covers the bundle of 
oriented tangent frames $\PSOM$ of $M$ such that the canonical diagram
$$
\xymatrix{
P \times \Spin(n) \ar[d] \ar[r]        & P \ar[d] \ar[dr] & \\
\PSOM \times \SO(n)  \ar[r]   & \PSOM \ar[r] & M
}
$$
commutes. 
Such a spin structure need not exist, e.g.\ complex projective plane 
$\CC\PP^2$ has none.
If $M$ has a spin structure we call $M$ a {\em spin manifold}.
The spin structure of a spin manifold is in general not unique.
More precisely, the cohomology $H^1(M;\ZZ_2)$ of a spin manifold acts simply 
transitively on the set of all spin structures.

Given a spin structure $P$ one can use the spinor representation
$$
\Spin(n) \to \Aut(\Sigma_n)
$$
to construct the associated {\em spinor bundle} $\Sigma M$ over $M$.
Here $\Sigma_n$ is a Hermitian vector space of dimension $2^{[n/2]}$ on which
$\Spin(n)$ acts by unitary transformations.
Hence $\Sigma M$ is a Hermitian vector bundle of rank $2^{[n/2]}$.
Sections in $\Sigma M$ are called {\em spinor fields} or simply {\em spinors}.
Note that unlike differential forms the definition of spinors requires the
choice of a spin structure.
The Levi-Civita connection on $\PSOM$ can be lifted to $P$ and therefore 
induces a covariant derivative $\nabla$ on $\Sigma M$.

Algebraic properties of the spinor representation ensure existence of
{\em Clifford multiplication}
$$
T_pM \otimes \Sigma_pM \to \Sigma_pM, 
\quad
X \otimes \psi \mapsto X \cdot \psi,
$$
satisfying the relations
$$
X\cdot Y \cdot \psi + Y \cdot X \cdot \psi + 2 \< X,Y\> \psi = 0
$$
for all $X,Y \in T_pM$, $\psi \in \Sigma_pM$, $p \in M$.
Here $\<\cdot,\cdot\>$ denotes the Riemannian metric.

The {\em Dirac operator} acting on spinors is defined as the composition
of $\nabla$ with Clifford multiplication.
Equivalently, if $e_1, \ldots,e_n$ is an orthonormal basis of $T_pM$, then
$$
(D\psi)(p) = \sum_{i=1}^n e_i\cdot\nabla_{e_i}\psi.
$$
The Dirac operator is a formally self-adjoint elliptic differential 
operator of first order.
If the underlying Riemannian manifold $M$ is complete, then $D$,
defined on compactly supported smooth spinors, is 
essentially self-adjoint in the Hilbert space of square-integrable
spinors. 
General elliptic theory ensures that the spectrum of $D$ is discrete
if $M$ is compact and satisfies Weyl's asymptotic law
$$
\lim\limits_{\lambda\rightarrow\infty}
\frac{N(\lambda)}{\lambda^n}
=
\frac{2^{[n/2]}\cdot\vol(M)}{(4\pi)^{\frac{n}{2}}
\cdot\Gamma\left(\frac{n}{2}+1\right)}.
$$
where $N(\lambda)$ is the number of eigenvalues whose modulus is $\le \lambda$.
This implies that the series
$$
\eta(s) = \sum_{\lambda\not= 0}\sign(\lambda)|\lambda|^{-s}
$$
converges for $s\in\CC$ if the real part of $s$ is sufficiently large.
Here summation is taken over all nonzero eigenvalues $\lambda$ of $D$,
each eigenvalue being repeated according to its multiplicity.
It can be shown that the function $\eta(s)$ extends to a meromorphic function
on the whole complex plane and has no pole at $s=0$.
Evaluation of this meromorphic extension at $s=0$ gives the
{$\eta$-invariant},
$$
\eta := \eta(0).
$$
If $M$ is complete but noncompact, then $D$ may also have eigenvalues of
infinite multiplicity, cumulation points of eigenvalues, and continuous
spectrum.

\section{The baby example}

In order to demonstrate the dependence of the Dirac spectrum
on the choice of spin structure the circle $S^1=\RR/2\pi\ZZ$ can serve as a
simple but nonetheless illustrative example.
Since the frame bundle $P_{\scriptsize \SO}S^1$ is trivial we can immediately
write down the trivial spin structure $P=S^1 \times \Spin(1)$.
Note that $\Spin(1)=\ZZ_2$ and $\Sigma_1=\CC$.
The associated spinor bundle is then also trivial and 1-dimensional.
Hence spinors are simply $\CC$-valued functions on $S^1$.
The Dirac operator is nothing but
$$
D = i\frac{d}{dt}.
$$
Elementary Fourier analysis shows that the spectrum consists of the
eigenvalues 
$$
\lambda_k = k
$$ 
with corresponding eigenfunction $t \mapsto
e^{-ikt}$, $k\in\ZZ$.
Since the spectrum is symmetric about zero, the $\eta$-series, and
in particular, the $\eta$-invariant vanishes,
$$
\eta=0.
$$
From $H^1(S^1;\ZZ_2)=\ZZ_2$ we see that $S^1$ has a second spin structure.
It can be described as $\tilde{P}=([0,2\pi]\times\Spin(1))/\sim$ where 
$\sim$ identifies $0$ with $2\pi$ while it interchanges the two
elements of $\Spin(1)$.
Let us call this spin structure the {\em nontrivial spin structure} of $S^1$.
Spinors with respect to this spin structure no longer correspond to 
functions on $S^1$, i.e.\ to $2\pi$-periodic functions on $\RR$, but rather
to $2\pi$-anti-periodic complex-valued functions on $\RR$,
$$
\psi(t+2\pi)=-\psi(t).
$$
This time the eigenvalues are 
$$
\lambda_k=k+\frac{1}{2},
$$
$k\in \ZZ$, with eigenfunction $t \mapsto e^{-i(k+\frac{1}{2})t}$.
Again, the spectrum is symmetric about $0$, hence $\eta=0$.
Vanishing of the $\eta$-invariant is in fact not surprising.
One can show that always $\eta=0$ for an $n$-dimensional manifold unless 
$n\equiv 3$ mod $4$.

The example $S^1$ has shown that the eigenvalues of the Dirac operator
definitely do depend on the choice of spin structure.
Even the dimension of the kernel of the Dirac operator is affected by a 
change of spin structure.
For the trivial spin structure of $S^1$ it is 1 while it is zero for the 
nontrivial spin structure.

We conclude this section with a remark on extendability of spin structures
because this sometimes causes confusion.
If $M$ is a Riemannian spin manifold with boundary $\partial M$, then
a spin structure on $M$ induces one on $\partial M$.
To see this consider the frame bundle $P_{\scriptsize \SO}\partial M$
of the boundary as a subbundle of $\PSOM$ restricted to the boundary
by completing a frame for $\partial M$ with the exterior unit normal vector
to a frame for $M$.
Now the inverse image of $P_{\scriptsize \SO}\partial M$ under the covering
map $P \to \PSOM$ defines a spin structure on $\partial M$.

Look at the case that $M$ is the disc with $S^1$ as its boundary.
Since the disk is simply connected it can have only one spin structure.
Hence only one of the two spin structures of $S^1$ extends to the disc.
The tangent vector to the boundary $S^1$ together with the unit normal vector
forms a frame for the disk which makes one full rotation when going around
the boundary one time.
It is therefore a loop in the frame bundle of the disk whose lift to the 
spin structure does not close up.
Thus the induced spin structure on the boundary is the nontrivial spin
structure of $S^1$ while the trivial spin structure does not bound.
Hence from a cobordism theoretical point of view the trivial spin structure
is nontrivial and vice versa.

\section{Flat tori and Bieberbach manifolds}

The case of higher-dimensional flat tori is very similar to the 1-dimensional
case.
There are $2^n$ different spin structures on $T^n=\RR^n/\Gamma$ where $\Gamma$
is a lattice in $\RR^n$.
Let $b_1, \ldots ,b_n$ be a basis of $\Gamma$, let $b_1^\ast, \ldots , 
b_n^\ast$ be the dual basis for the dual lattice $\Gamma^\ast$.
Spin structures can then be classified by $n$-tuples $(\delta_1, \ldots, 
\delta_n)$ where each $\delta_j \in \{ 0,1\}$ indicates whether or not the
spin structure is twisted in direction $b_j$.
The spectrum of the Dirac operator can then be computed:

\begin{theo}[Friedrich \cite{friedrich84a}]
\label{torus}
The eigenvalues of the Dirac operator on $T^n=\RR^n/\Gamma$ with spin structure
corresponding to $(\delta_1, \ldots, \delta_n)$ are given by
$$
\pm 2\pi \left| b^\ast + \frac{1}{2}\sum_{j=1}^n \delta_j b_j^\ast  \right|
$$
where $b^\ast$ runs through $\Gamma^\ast$ and each $b^\ast$ contributes 
multiplicity $2^{[n/2]-1}$.
\end{theo}

Again the spectrum depends on the choice of spin structure.
In particular, eigenvalue $0$ occurs only for the trivial spin structure
given by $(\delta_1, \ldots, \delta_n) = (0,\ldots,0)$.
Since again the spectrum is symmetric about zero, the $\eta$-invariant
vanishes, $\eta=0$, for all spin structures.

This changes if one passes from tori to more general compact connected flat 
manifolds, also called {\em Bieberbach manifolds}.
They can always be written as a quotient $M = G\backslash T^n$ of a torus by 
a finite group $G$.
In three dimensions, $n=3$, there are 5 classes of compact oriented Bieberbach
manifolds besides the torus.
Their Dirac spectra have been calculated by Pf\"affle \cite{pfaeffle99ppa}
for all flat metrics.
This time one finds examples with asymmetric spectrum and the $\eta$-invariant
depends on the choice of spin structure.

\begin{theo}[Pf\"affle \cite{pfaeffle99ppa}]
The $\eta$-invariant of the $3$-dimensional compact oriented Bieberbach
manifolds besides the torus are given by the following table:
\begin{center}
\begin{tabular}{|c|c|c|c|c|}
\hline
$G$ & total \# spin structures & \multicolumn{3}{c|}{$\eta$-invariant for \#
spin structures}\\
\hline\hline
$\ZZ_2$ & 8 & $\eta=0$ for 6 & $\eta=1$ for 1 & $\eta=-1$ for 1 \\
\hline
$\ZZ_3$ & 2 &  & $\eta=\frac{4}{3}$ for 1 & $\eta=-\frac{2}{3}$ for 1 \\
\hline
$\ZZ_4$ & 4 & $\eta=0$ for 2 & $\eta=\frac{3}{2}$ for 1 & $\eta=-\frac{1}{2}$ for 1 \\
\hline
$\ZZ_6$ & 2 &  & $\eta=\frac{5}{3}$ for 1 & $\eta=-\frac{1}{3}$ for 1 \\
\hline
$\ZZ_2 \times \ZZ_2$ & 4 & $\eta=0$ for 4 & & \\
\hline
\end{tabular}
\end{center}
\end{theo}
\begin{center}
{\rm Table 1}
\end{center}

Note that the $\eta$-invariant does not depend on the choice
of flat metric even though the spectrum does.
Depending on $G$ there is a 2-, 3- or 4-parameter family of such metrics on 
$M$.

\section{Spherical space forms}

The Dirac spectrum on the sphere $S^n$ with constant curvature has been
computed by different methods in \cite{baer96a,sulanke79a,trautman95a}.
The eigenvalues are
\begin{equation}
\pm\left(\frac{n}{2}+k \right),
\label{sneigenwert}
\end{equation}
$k\in\NN_0$, with multiplicity $2^{[n/2]}\cdot\left(\begin{array}{c}k+n-1 \\
k\end{array}\right)$.
For $n\ge 2$ the sphere is simply connected, hence has only one spin structure.
Therefore let us look at spherical space forms $M=\Gamma\backslash S^n$ where 
$\Gamma$ is a finite fixed point free subgroup of $\SO(n+1)$.
Spin structures correspond to homomorphisms $\epsilon : \Gamma \to \Spin(n+1)$
such that
$$
\xymatrix{
           & \Spin(n+1) \ar[d] \\
\Gamma \ar[ur]^\epsilon \ar[r]   & \SO(n+1)
}
$$
commutes.
Since any eigenspinor on $M$ can be lifted to $S^n$ all eigenvalues of $M$
are also eigenvalues of $S^n$, hence of the form (\ref{sneigenwert}).
To know the spectrum of $M$ one must compute the multiplicities $\mu_k$ of
$\frac{n}{2}+k$ and $\mu_{-k}$ of $-(\frac{n}{2}+k)$.
They can be most easily expressed by encoding them into two power series,
so-called {\em Poincar\'e series}
\begin{eqnarray*}
F_+(z) &=& \sum_{k=0}^\infty \mu_k z^k,\\
F_-(z) &=& \sum_{k=0}^\infty \mu_{-k} z^k.
\end{eqnarray*}
To formulate the result recall that in even dimension $2m$ the spinor 
representation is reducible and can be decomposed into two half spinor
representations
$$
\Spin(n) \to \Aut(\Sigma_{2m}^\pm),
$$
$\Sigma_{2m}= \Sigma_{2m}^+ \oplus \Sigma_{2m}^-$.
Denote their characters by $\chi^\pm : \Spin(2m) \to \CC$.

\begin{theo}[B\"ar \cite{baer96a}]
\label{ssfeigenwerte}
Let $M=\Gamma \backslash S^{n}$, $n=2m-1$, be a spherical space form with
spin structure given by $\epsilon : \Gamma \to \Spin(2m)$.
Then the eigenvalues of the Dirac operator are $\pm (\frac{n}{2} + k), k
\geq 0,$ with multiplicities determined by
\begin{eqnarray*}
F_+(z) & = & \frac{1}{|\Gamma |} \sum_{\gamma \in \Gamma} \frac{\chi^-
(\epsilon (\gamma )) - z \cdot \chi^+(\epsilon (\gamma ))}{\det(1_{2m} -
z\cdot\gamma )} , \\
F_-(z) & = & \frac{1}{|\Gamma |} \sum_{\gamma \in \Gamma} \frac{\chi^+
(\epsilon (\gamma )) - z \cdot \chi^-(\epsilon (\gamma ))}{\det(1_{2m} -
z\cdot\gamma )} .
\end{eqnarray*}
\end{theo}

Note that only odd-dimensional spherical space forms are of interest because
in even dimensions real projective space is the only quotient and in this case
it is not even orientable.

Let us use Theorem \ref{ssfeigenwerte} to compute the $\eta$-invariant of
spherical space forms.
We get immediately for the {\em $\theta$-functions}
\begin{eqnarray*}
\theta_{\pm} (t) &:=& e^{- \frac{n}{2}t} \cdot F_{\pm} (e^{-t}) \\
&=&\frac{e^{-(m+\frac{1}{2})t}}{|\Gamma|} \sum\limits_{\gamma\in\Gamma} 
\frac{\chi^\mp (\epsilon(\gamma)) - e^{-t} \cdot \chi^{\pm} (\epsilon(\gamma))}
{\det (1_{2m} - e^{-t} \cdot \gamma)}.
\end{eqnarray*}
The coefficient of $t^0$ in the Laurent expansion at $t=0$ is given by 
\begin{eqnarray*}
LR_0 (\theta_+)  
&=& 
\frac{1}{|\Gamma|} 
\sum\limits_{\gamma\in\Gamma - \left\{1_{2m}\right\}} 
\frac{\chi^- (\epsilon(\gamma)) - \chi^+ (\epsilon(\gamma))}
{\det (1_{2m}-\gamma)} \\
&& + \hspace{0.2cm}
LR_0 \left(
\frac{e^{-(m+\frac{1}{2})t}}{|\Gamma|}
\cdot
\frac{2^{m-1} -e^{-t} \cdot 2^{m-1}}{\det(1_{2m}-e^{-t} \cdot 1_{2m})}
\right).
\end{eqnarray*}
Similarly, 
\begin{eqnarray*}
LR_0 (\theta_-)  
&=& 
\frac{1}{|\Gamma|} 
\sum\limits_{\gamma\in\Gamma - \left\{1_{2m}\right\}} 
\frac{\chi^+ (\epsilon(\gamma)) - \chi^- (\epsilon(\gamma))}
{\det (1_{2m}-\gamma)} \\
&&  +  \hspace{0.2cm}
LR_0 \left(
\frac{e^{-(m+\frac{1}{2})t}}{|\Gamma|}
\cdot
\frac{2^{m-1} -e^{-t} \cdot 2^{m-1}}{\det(1_{2m}-e^{-t} \cdot 1_{2m})}
\right).
\end{eqnarray*}
Hence we obtain for $\theta := \theta_+ - \theta_-$
\begin{eqnarray*}
LR_0 (\theta)  
&=& 
LR_0 (\theta_+) - LR_0 (\theta_-) \\
&=&
\frac{2}{|\Gamma|} \sum\limits_{\gamma\in\Gamma - \left\{1_{2m}\right\}} 
\frac{\chi^- (\epsilon(\gamma)) - \chi^+ (\epsilon(\gamma))}
{\det (1_{2m}-\gamma)}.
\end{eqnarray*} 
The same argument shows that the poles of $\theta_+$ and $\theta_-$ cancel,
hence $\theta$ is holomorphic at $t=0$ with
$$
\theta(0) = 
\frac{2}{|\Gamma|} \sum\limits_{\gamma\in\Gamma - \left\{1_{2m}\right\}} 
\frac{(\chi^- - \chi^+)(\epsilon(\gamma))}{\det (1_{2m}-\gamma)}.
$$
Now we observe that
$$
\theta_+(t) = \sum_{k=0}^\infty \mu_k e^{-(n/2 + k)t}
= \sum_{\lambda > 0} e^{-\lambda t},
$$
and similarly for $\theta_-$.
Application of the {\em Mellin transformation} yields
$$
\eta(s) = \frac{1}{\Gamma(s)}\int_0^\infty \theta(t) t^{s-1} dt.
$$
Therefore 
$$
\eta = \lim_{s\to 0}\frac{1}{\Gamma(s)}\int_0^\infty \theta(t) t^{s-1} dt
= \Res_{s=0}\left(\int_0^\infty \theta(t) t^{s-1} dt\right).
$$
Since $\theta$ decays exponentially fast for $t \to \infty$ the function
$s \mapsto \int_1^\infty \theta(t) t^{s-1} dt$ is holomorphic at $s=0$.
Thus 
$$
\eta = \Res_{s=0}\left(\int_0^1 \theta(t) t^{s-1} dt\right) = \theta(0).
$$
We have proved 

\begin{theo}
Let $M=\Gamma \backslash S^{2m-1}$ be a spherical space form with
spin structure given by $\epsilon : \Gamma \to \Spin(2m)$.
Then the $\eta$-invariant of $M$ is given by
$$
\eta = \frac{2}{|\Gamma|} 
\sum\limits_{\gamma\in\Gamma-\left\{1_{2m}\right\}} 
\frac{(\chi^- - \chi^+)(\epsilon(\gamma))}
{\det (1_{2m}-\gamma)}.
$$
\end{theo}
 
\begin{exem}
We take a look at real projective space $\RR\PP^{2m-1}$, i.e. 
$\Gamma=\left\{1_{2m}, -1_{2m}\right\}$. 
If we view $\Spin(2m)$ as sitting in the Clifford algebra 
$\Cl(\RR^{2m})$, compare \cite{lawson-michelsohn89a}, then we can define the 
``volume element'' 
$$
\omega:=e_1\cdot e_2\cdot\ldots\cdot e_{2m}\in\Spin(2m)\subset\Cl(\RR^{2m}) 
$$
where $e_1, \ldots, e_{2m}$ denotes the standard basis of $\RR^{2m}$. 
It is not hard to see that under the map $\Spin(2m) \to \SO(2m)$ the
volume element $\omega$ is mapped to $-1_{2m}$. 
Hence the two preimages of $-1_{2m}$ in $\Spin(2m)$ are $\pm\omega$. 
To specify a spin structure we may define 
$$
\epsilon(-1_{2m}):=\omega
$$
or
$$
\epsilon(-1_{2m}):=-\omega.
$$
One checks
$$
1=\epsilon(1_{2m}) = \epsilon((-1_{2m})^2) = \epsilon(-1_{2m})^2 = 
(\pm\omega)^2 = (-1)^m.
$$ 
Hence $\RR\PP^{2m-1}$ is not spin if $m$ is odd $(m\ge 3)$ whereas it has 
two spin structures if $m$ is even. 

The volume element $\omega$ acts on the two half spinor spaces via 
multiplication by $\pm 1$. 
Hence
\begin{eqnarray*}
\eta &=& \pm\frac{2}{2}\cdot\frac{2\cdot 2^{m-1}}{2^{2m}} \\
     &=& \pm 2^{-m}
\end{eqnarray*} 
We summarize 
\end{exem}

\begin{coro}
\label{rpn} 
For $n\ge 2$ real projective space $\RR\PP^n$ is spin if and only if
$n\equiv 3$ mod 4, in which case it has exactly two spin structures. 
The $\eta$-invariant for the Dirac operator is given by
$$
\eta = \pm 2^{-m}, n=2m-1,
$$
where the sign depends on the spin structure chosen.
\end{coro}

See also \cite{gilkey84a,gilkey89a} where the $\eta$-invariant of all
twisted signature operators on spherical space forms is determined and used 
to compute their $K$-theory. 

\section{Eigenvalue estimates}

Up to very recently all known lower eigenvalue estimates for the
Dirac operator did not take into account the spin structure despite
its influence on the spectrum that we have encountered in the
examples.
This is due to the fact that they are all based on the Bochner technique,
hence on a {\em local} computation.
To find estimates which can see the spin structure one needs to define
new, truly spin geometric invariants.
Such invariants have been proposed by Ammann \cite{ammann98a,ammann00ppa} in 
the case of a 2-torus.
Recall the definition of the {\em systole} of a Riemannian manifold $(M,g)$
$$
\sys_1(M,g) = \inf\{ \length(\gamma)\ |\ \gamma \mbox{ is a noncontractible 
loop}  \}.
$$
In case $M$ is a torus there is a canonical spin structure, the trivial
spin structure $P_0$.
Hence the set of spin structures can be identified with $H^1(M;\ZZ_2)$
by identifying $P_0$ with $0$.
It then makes sense to evaluate a spin structure $P$ on the homology
class of a loop $\gamma$ yielding an element in $\ZZ_2 = \{ 1,-1 \}$.
This value $P([\gamma])$ specifies whether or not the spin structure $P$
is twisted along $\gamma$.
Ammann defines the {\em spinning systole}
$$
\spinsys_1(M,g,P) = \inf\{ \length(\gamma)\ |\ \gamma \mbox{ is a loop with }
P([\gamma]) = -1 \}.
$$
Hence the infimum is taken only over those loops along which the spin structure
twists.
In case the spin structure is trivial, $P=P_0$, the spinning systole is
infinite.

\begin{theo}[Ammann \cite{ammann00ppa}]
Let $g$ be a Riemannian metric on the 2-torus whose Gauss curvature $K$
satisfies $\| K\|_{L^1(T^2,g)}<4\pi$.
Let $P$ be a spin structure on $T^2$.

Then for all eigenvalues $\lambda$ of the Dirac operator the estimate
$$
\lambda^2 \ge \frac{C(\| K\|_{L^1(T^2,g)},\| K\|_{L^2(T^2,g)},\area(T^2,g),
\sys_1(T^2,g))}{\spinsys_1(T^2,g,P)^2}
$$
holds where $C(\| K\|_{L^1(T^2,g)},\| K\|_{L^2(T^2,g)},\area(T^2,g),
\sys_1(T^2,g))>0$ is an explicitly given expression.
\end{theo}

The estimate is sharp in the sense that for some flat metrics equality is
attained.
In a similar way Ammann defines the {\em nonspinning systole} and proves
an analogous estimate for which equality is attaind for all flat metrics.
The proofs are based on a comparison of the Dirac spectra for the metric
$g$ with the one for the conformally equivalent flat metric $g_0$.
Remember that by Theorem \ref{torus} the spectrum for $g_0$ is 
explicitly known.
Most of the work is then done to control the oscillation of the function
which relates the two conformally equivalent metrics $g$ and $g_0$ in terms
of the geometric data occuring in $C(\| K\|_{L^1(T^2,g)},\| K\|_{L^2(T^2,g)},
\area(T^2,g),\sys_1(T^2,g))$.
This way it is also possible to derive upper eigenvalue estimates, see
\cite{ammann98a,ammann00ppa} for details.

In the same paper \cite{ammann00ppa} Ammann also studies the question
how far Dirac spectra for different metrics on a compact manifold can
be away from each other.
If $P_1$ and $P_2$ are two spin structures on a Riemannian manifold $(M,g)$,
then there is a unique $\chi\in H^1(M;\ZZ_2)$ taking $P_1$ to $P_2$.
On $H^1(M;\ZZ_2)$ there is a canonical norm, the {\em stable norm} (or
$L^\infty$-norm).
Ammann shows that if the Dirac eigenvalues $\lambda_j$ of $(M,g,P_1)$
and $\lambda'_j$ of $(M,g,P_2)$ are numbered correctly, then
$$
|\lambda_j - \lambda'_j| \le 2\pi\|\chi\|_{L^\infty}.
$$

\section{Collapse of circle bundles}

Another instance where the choice of spin structure has strong influence
on the spectral behavior occurs when one looks at circle bundles and
their collapse to the basis.
To this extent let $(M,g_M)$ be a compact Riemannian spin manifold with an 
isometric and free circle action.
For simplicity we suppose that the fibers have constant lengths.
We give the quotient $N:=S^1\backslash M$ the unique Riemannian metric $g_N$
for which the projection $M \to N$ is a Riemannian submersion.
By rescaling the metric $g_M$ along the fibers while keeping it unchanged
on the orthogonal complement to the fibers we obtain a 1-parameter family
of Riemannian metrics $g_\ell$ on $M$ with respect to which $(M,g_\ell)\to 
(N,g_N)$ is a Riemannian submersion and the fibers are of length $2\pi\ell$.
{\em Collapse} of this circle bundle now means that we let $\ell\to 0$, i.e.\
we shrink the fibers to a point.
Then $(M,g_\ell)$ tends to $(N,g_N)$ in the Gromov-Hausdorff topology.
In the physics literature this is also refered to as {\em adiabatic limit}.
The question now is how the spectrum behaves.
In particular, do eigenvalues of $(M,g_\ell)$ tend to those of $(N,g_N)$?

For the answer we have to study the spin structure $P$ on $M$.
The isometric circle action on $M$ induces a circle action on the frame bundle
$\PSOM$.
This $S^1$-action may or may not lift to $P$.
In case it lifts we call $P$ {\em projectable}, otherwise we call it 
{\em nonprojectable}.
If $P$ is projectable, then it induces a spin structure on $N$.

\begin{theo}[Ammann-B\"ar \cite{ammann-baer98a}]
\label{projectable}
Let $P$ be projectable and let $N$ carry the induced spin structure.
Denote the Dirac eigenvalues of $(N,g_N)$ by $\mu_j$.
Then the Dirac eigenvalues $\lambda_{j,k}(\ell)$, $j,k\in\ZZ$, of $(M,g_\ell)$,
if numbered correctly, depend continuously on $\ell$ and for
$\ell \to 0$ the following holds:
\begin{itemize}
\item
For all $j$ and $k$ we have
$$
\ell\cdot\lambda_{j,k}(\ell) \to k.
$$
In particular, $\lambda_{j,k}(\ell) \to \pm\infty$ for $k\not= 0$.
\item
If $\dim(N)$ is even, then
$$
\lambda_{j,0}(\ell) \to \mu_j.
$$
\item
If $\dim(N)$ is odd, then
$$
\lambda_{2j-1,0}(\ell) \to \mu_j,
\lambda_{2j,0}(\ell) \to -\mu_j.
$$
\end{itemize}
\end{theo}

Roughly, some eigenvalues tend to $\pm\infty$ while the others converge
to the eigenvalues of the bases (and their negatives for odd-dimensional basis).
This can be applied to the Hopf fibration $S^{2m+1}\to \CC\PP^m$.
If $m$ is odd, then the unique spin structure on $S^{2m+1}$ is projectable
and one can use Theorem \ref{projectable} to compute the spectrum of 
complex projective space. 
If $m$ is even, then the spin structure on $S^{2m+1}$ is not projectable.
Indeed $\CC\PP^m$ is not spin in this case.
The behavior of the spectrum is in this case described by the following

\begin{theo}[Ammann-B\"ar \cite{ammann-baer98a}]
Let $P$ be nonprojectable.
Then the Dirac eigenvalues $\lambda_{j,k}(\ell)$, $j\in\ZZ$, $k\in\ZZ+(1/2)$, 
of $(M,g_\ell)$, if numbered correctly, depend continuously 
on $\ell$ and for $\ell \to 0$ the following holds:
For all $j$ and $k$ we have
$$
\ell\cdot\lambda_{j,k}(\ell) \to k.
$$
In particular, $\lambda_{j,k}(\ell) \to \pm\infty$ for all $k$ and $j$.
\end{theo}

Both cases occur e.g.\ for Heisenberg manifolds.
They are circle bundles over flat tori.
The proofs are based on a Fourier decomposition along the fibers.
For the case varying fiber length see \cite{ammann99a}, for a very recent
paper containing a quite general treatment of collapse see \cite{lott00ppa}.

\section{$\eta$-invariant}

We have already seen in examples that the $\eta$-invariant does depend
on the spin structure.
However it turns out that the $\eta$-invariants for different spin
structures on the same Riemannian manifold $M$ are not totally unrelated.
Recall that for two spin structures $P_1$ and $P_2$ there is a unique
$\chi\in H^1(M;\ZZ_2)$ mapping $P_1$ to $P_2$.
We call $\chi$ {\em realizable as a differential form} if there exists
a 1-form $\omega$ such that 
$$
\exp\left(2\pi i\int_\gamma\omega\right) = \chi([\gamma])
$$
for all loops $\gamma$.
This is equivalent to the vanishing of $\chi$ on the mod-2-reduction of
all torsion elements in $H^1(M;\ZZ)$.
See \cite{ammann98a} for this and other characterizations.

\begin{theo}[Dahl \cite{dahl99ppa}]
\label{etadiff}
Let $P_1$ and $P_2$ be two spin structures on the compact Riemannian
manifold $M$.
Suppose the element $\chi\in H^1(M;\ZZ_2)$ mapping $P_1$ to $P_2$ is 
realizable as a differential form.
Then
$$
\eta_{M,P_1} - \eta_{M,P_2} \in \ZZ.
$$
\end{theo}

Be careful that some conventions and in particular the definition of the
$\eta$-invariant in \cite{dahl99ppa} differ from ours.
One can check that in the case of 3-dimensional Bieberbach manifolds
the assumption on $\chi\in H^1(G\backslash T^3;\ZZ_2)$ is always fulfilled for 
$G=\ZZ_3$ and for $G=\ZZ_2\times\ZZ_2$.
It is fulfilled for some but not all $\chi\in H^1(G\backslash T^3;\ZZ_2)$ in 
case $G=\ZZ_2$ and $G=\ZZ_4$.

From $H_1(\RR\PP^n;\ZZ)=\ZZ_2$ one sees that the nontrivial element of
$H^1(\RR\PP^n;\ZZ_2)$ is not realizable as a differential form.
In fact, otherwise Theorem \ref{etadiff} would contradict Corollary \ref{rpn}.
This example shows that this assumption on $\chi$ cannot be dispensed with.

The proof of Dahl's theorem is based on a suitable application of the
Atiyah-Patodi-Singer index theorem \cite{atiyah-patodi-singer75a} to the
cylinder over $M$.
The main idea is to write the difference of $\eta$-invariants as a 
linear combination of indices, hence of integers.
This index theorem was the reason to introduce the $\eta$-invariant 
in the first place.

\section{Noncompact hyperbolic manifolds}

In contrast to spaces of constant sectional curvature $\ge 0$ there
is no hope to be able to explicitly compute the Dirac spectrum
on a space of constant negative curvature.
In \cite{baer-schmutz92a,bures93a,bures98a,hitchin74a} the dimension of the 
kernel of the Dirac operator on hyperbolic Riemann surfaces is considered.
For hyperelliptic metrics it can be computed for all spin structures
and it varies with the spin structure. 

So far we only have considered {\em compact} manifolds whose Dirac spectrum
is always discrete.
Let us now discuss noncompact hyperbolic manifolds with an eye to the
question whether or not the continuous spectrum also depends on the
choice of spin structure.

A {\em hyperbolic manifold} is a complete connected Riemannian
manifold of constant sectional curvature -1.
Every hyperbolic manifold $M$ of finite volume can be decomposed disjointly 
into a relatively compact $M_0$ and finitely many cusps $\E_j$,
$$
M = M_0 \sqcup \bigsqcup_{j=1}^k \E_j
$$
\begin{center}
\pspicture(1,0)(14,6)

\pscustom[linecolor=white,fillstyle=solid,fillcolor=red]{
  \pscurve(4.95,6)(4.9,5.5)(4.7,5)(4.4,4.5)
  \pscurve(4.4,4.5)(4.6,4.4)(5.05,4.3)(5.7,4.5)
  \pscurve(5.7,4.5)(5.4,5)(5.2,5.5)(5.15,6)
  \pscurve(5.15,6)(4.95,6)}
\pscustom[linecolor=white,fillstyle=solid,fillcolor=red]{
  \pscurve(7,3.6)(8,3.35)(9,3.22)(10,3.16)(11,3.13)
  \pscurve(11,3.13)(11,2.97)
  \pscurve(11,2.96)(10.6,2.96)(10,2.94)(9,2.88)(8,2.75)(7,2.5)
  \pscurve(7,2.5)(6.9,2.6)(6.8,3.05)(7,3.6)}

\pscurve(7,2.5)(8,2.75)(9,2.88)(10,2.94)(11,2.97)
\pscurve(7,3.6)(8,3.35)(9,3.22)(10,3.16)(11,3.13)
\pscurve(4.4,4.5)(4.7,5)(4.9,5.5)(4.95,6)
\pscurve(5.7,4.5)(5.4,5)(5.2,5.5)(5.15,6)
\psecurve(5.4,5)(5.7,4.5)(6.15,3.95)(7,3.6)(8,3.35)
\psecurve(4.7,5)(4.4,4.5)(3,3)(3,1)(4.5,0.7)(5.9,1.2)(7,2.5)(8,2.75)
\pscurve(4.1,2.75)(4,2.5)(4.5,2.1)(5,2.2)(5.2,2.3)
\pscurve(4,2.5)(4.6,2.45)(5,2.2)
\pscurve(7,2.5)(6.8,3.05)(7,3.6)
\pscurve(4.4,4.5)(5.05,4.3)(5.7,4.5)

\rput(5.3,2.9){$M_0$}
\rput(7.3,3.1){\psframebox*[framearc=0.5]{$\E_1$}}
\rput(5.05,4.8){\psframebox*[framearc=0.5]{$\E_2$}}
\rput(6.55,3.1){{$N_1$}}
\rput(5.05,4.05){{$N_2$}}

\endpspicture
\vspace{-1cm}
\end{center}
\begin{center}
\rm Fig.~1
\end{center}
\noindent
where each $\E_j$ is of the form $\E_j = N_j \times [0,\infty)$.
Here $N_j$ denotes a connected compact manifold with a flat metric 
$g_{N_j}$, a Bieberbach manifold, and $\E_j$ carries the warped 
product metric $g_{\E_j} = e^{-2t}\cdot g_{N_j} + dt^2$.
If $M$ is 2- or 3-dimensional and oriented, then $N_j$ is a circle $S^1$ or 
a 2-torus $T^2$ respectively.
We call a spin structure on $M$ {\em trivial along the cusp $\E_j$} if 
its restriction to $N_j$ yields the trivial spin structure on $N_j$.
Otherwise we call it {\em nontrivial along $\E_j$}.

Now it turns out that only two extremal cases occur for the spectrum of
the Dirac operator, it is either discrete as in the compact case or it is
the whole real line.
And it is the spin structure which is responsible for the choice between the
two cases.

\begin{theo}[B\"ar \cite{baer98ppb}]
\label{dichotomy}
Let $M$ be a hyperbolic 2- or 3-manifold of finite volume equipped with
a spin structure.

If the spin structure is trivial along at least one cusp, then the Dirac 
spectrum is the whole real line
$$
spec(D) = \RR .
$$

If the spin structure is nontrivial along all cusps, then the spectrum
is discrete.
\end{theo}

In fact, this theorem also holds in higher dimensions.
The proof is based on the fact that the essential spectrum of the Dirac
operator is unaffected by changes in compact regions.
Hence one only needs to look at the cusps and they are given in a very 
explicit form.
A separation of variables along the cusps yields the result.
Of course, Theorem \ref{dichotomy} does not say anything about existence
of spin structures on $M$ being trivial or nontrivial along the various cusps.
This can be examined by topological methods and the answer for hyperbolic
surfaces is given in the table
\begin{center}
\begin{tabular}{|c||c|c|}
\hline
\multicolumn{3}{|c|}{Hyperbolic surface of finite area}\\
\hline\hline
\# of cusps & 
\begin{tabular}{c}
existence of spin structure\\ 
with discrete spectrum 
\end{tabular}
&
\begin{tabular}{c}
existence of spin structure\\ 
with $spec(D) = \RR$
\end{tabular} \\
\hline\hline
0 & YES & NO\\
\hline
1 & YES & NO\\
\hline
$\ge 2$ & YES & YES\\
\hline
\end{tabular}
\end{center}
\begin{center}
{\rm Table 2}
\end{center}

\noindent
while the 3-dimensional case is given by

\begin{center}
\begin{tabular}{|c||c|c|}
\hline
\multicolumn{3}{|c|}{Hyperbolic 3-manifold of finite volume}\\
\hline\hline
\# of cusps & 
\begin{tabular}{c}
existence of spin structure\\ 
with discrete spectrum 
\end{tabular}
&
\begin{tabular}{c}
existence of spin structure\\ 
with $spec(D) = \RR$
\end{tabular} \\
\hline\hline
0 & YES & NO\\
\hline
1 & YES & NO\\
\hline
$\ge 2$ & YES & depends on $M$\\
\hline
\end{tabular}
\end{center}
\begin{center}
{\rm Table 3}
\end{center}

The tables show that hyperbolic 2- or 3-manifolds of finite volume with
one end behave like compact ones, the Dirac spectrum is always discrete.
A surface with two or more ends always admits both types of spin structures.
This is not true for 3-manifolds.
Discrete spectrum is always possible but the case $spec(D) = \RR$ only
sometimes.
If the hyperbolic 3-manifold is topologically given as the complement
of a link in $S^3$ (and this construction is one of the main sources for
hyperbolic 3-manifolds of finite volume), then this question can be decided.

\begin{theo}[B\"ar \cite{baer98ppb}]
Let $K \subset S^3$ be a link, let $M = S^3 - K$ carry a hyperbolic metric
of finite volume.

If the linking number of all pairs of components $(K_i,K_j)$ of $K$
is even,
$$
Lk(K_i,K_j) \equiv 0 \mbox{ mod } 2,
$$
$i \not= j$, then the spectrum of the Dirac operator on $M$ is discrete for all
spin structures.

If there exist two components $K_i$ and $K_j$ of $K$, $i \not= j$, with
odd linking number, then $M$ has a spin structure such that the spectrum of the
Dirac operator satisfies
$$
spec(D) = \RR.
$$
\end{theo}

The condition on the linking numbers is very easy to
verify in given examples.
Since we compute modulo 2 orientations of link components are irrelevant.
If the link is given by a planar projection, then modulo 2, $Lk(K_i,K_j)$
is the same as the number of over-crossings of $K_i$ over $K_j$.

\begin{exem}
The complements of the following links possess a hyperbolic structure
of finite volume.
All linking numbers are even.
Hence the Dirac spectrum on those hyperbolic manifolds is discrete
for all spin structures.

\begin{center}
\pspicture(1,0)(14,10)

\psecurve[showpoints=false,linewidth=2pt,linecolor=blue](2.5,7)(2.6,6.6)(3,6.1)(3.5,5.9)
(4.1,6.2)(4.4,6.8)(4.5,7.4)(4.4,8)(4.3,8.4)
\psecurve[showpoints=false,linewidth=2pt,linecolor=blue](4.4,8)(4.3,8.4)(4.1,8.8)(3.6,9)
(3,8.8)(2.6,8.2)(2.5,7.6)(2.5,7)(2.6,6.6)

\psecurve[showpoints=false,linewidth=2pt,linecolor=red](3.3,7.3)(3.7,7.7)(3.9,8)(4.4,8.2)
(5.1,7.8)(5.1,7.2)(4.6,6.9)(4.2,6.9)
\psecurve[showpoints=false,linewidth=2pt,linecolor=red](4.6,6.9)(4.2,6.9)(3.8,7.2)(3.5,7.5)
(3.2,7.8)(2.7,8.1)(2.4,8.1)
\psecurve[showpoints=false,linewidth=2pt,linecolor=red](2.7,8.1)(2.4,8.1)(1.9,7.8)(1.9,7.2)
(2.5,6.8)(3,6.9)(3.3,7.3)(3.7,7.7)

\rput(5.5,9.5){$5^2_1$}


\psecurve[showpoints=false,linewidth=2pt,linecolor=blue](9.5,7)(9.6,6.6)(10,6.1)(10.5,5.9)
(11.1,6.2)(11.4,6.8)(11.5,7.4)(11.4,8)(11.3,8.4)
\psecurve[showpoints=false,linewidth=2pt,linecolor=blue](11.4,8)(11.3,8.4)(11.1,8.8)(10.6,9)
(10,8.8)(9.6,8.2)(9.5,7.6)(9.5,7)(9.6,6.6)

\psecurve[showpoints=false,linewidth=2pt,linecolor=red](11.3,7.1)(11.6,7)(12,7.4)(12,7.8)
(11.6,8.2)(11.1,8.1)(10.9,7.8)(10.8,7.4)
\psecurve[showpoints=false,linewidth=2pt,linecolor=red](10.9,7.8)(10.8,7.4)(10.5,7.1)
(10.2,7.2)(10,7.8)(9.7,8)(9.4,8)
\psecurve[showpoints=false,linewidth=2pt,linecolor=red](9.6,7.9)(9.4,8)(9,7.6)(9,7)
(9.4,6.75)(10,7.3)(10.2,7.6)
\psecurve[showpoints=false,linewidth=2pt,linecolor=red](10,7.3)(10.2,7.6)(10.5,7.8)
(10.9,7.6)(11.3,7.1)(11.6,7)

\rput(12.5,9.5){$6^2_3$}


\psecurve[showpoints=false,linewidth=2pt,linecolor=blue](2.8,3.8)(2.9,4.05)(3.4,4.3)
(3.9,4)(4,3.4)(3.9,3)(3.8,2.8)
\psecurve[showpoints=false,linewidth=2pt,linecolor=blue](3.9,3)(3.8,2.8)(3.6,2.6)
(3.2,2.6)(2.8,3)(2.7,3.4)(2.8,3.8)(2.9,4.05)

\psecurve[showpoints=false,linewidth=2pt,linecolor=red](4.1,3.9)(3.8,3.9)(3,3.9)
(2.4,3.8)(2,3.2)(2.4,2.4)(2.6,2.1)
\psecurve[showpoints=false,linewidth=2pt,linecolor=red](2.4,2.4)(2.6,2.1)(3,1.6)
(3.6,1.2)(4.4,1.4)(4.5,2)(4.3,2.2)
\psecurve[showpoints=false,linewidth=2pt,linecolor=red](4.5,2)(4.3,2.2)(3.8,3)
(3.4,3.1)(3,3)(2.8,2.7)
\psecurve[showpoints=false,linewidth=2pt,linecolor=red](3,3)(2.8,2.7)(2.4,2.2)
(2.2,1.6)(2.8,1.1)(3.3,1.2)(3.6,1.4)
\psecurve[showpoints=false,linewidth=2pt,linecolor=red](3.3,1.2)(3.6,1.4)(4.3,2)
(4.8,2.8)(4.8,3.4)(4.4,3.8)(4.1,3.9)(3.8,3.9)

\rput(5.5,4.5){$7^2_4$}


\psarc[linewidth=2pt,linecolor=blue](11,2){1}{70}{163}
\psarc[linewidth=2pt,linecolor=blue](11,2){1}{183}{50}

\psarc[linewidth=2pt,linecolor=red](9.8,2){1}{65}{296}
\psarc[linewidth=2pt,linecolor=red](9.8,2){1}{316}{45}

\psarc[linewidth=2pt,linecolor=green](10.5,3){1}{296}{173}
\psarc[linewidth=2pt,linecolor=green](10.5,3){1}{193}{276}

\rput(12.5,4.5){$6^3_2$}


\endpspicture
\end{center}

\begin{center}
\rm Fig.~2
\end{center}
\noindent
This example includes the Whitehead link ($5^2_1$) and the Borromeo rings
($6^3_2$).
\end{exem}

\begin{exem}
The complements of the following links possess a hyperbolic structure
of finite volume.
There are odd linking numbers.
Hence those hyperbolic manifolds have a spin structure for which the 
Dirac spectrum is the whole real line.

\begin{center}
\pspicture(1,0)(14,10)

\psecurve[showpoints=false,linewidth=2pt,linecolor=blue](2.6,7)(2.8,6.8)(3.4,6.4)
(3.8,6.2)(4.6,6.2)(4.7,6.8)(4.4,7)
\psecurve[showpoints=false,linewidth=2pt,linecolor=blue](4.7,6.8)(4.4,7)(4,7.2)
(3.6,7.4)(3.3,7.8)(3.4,8.2)(3.7,8.5)
\psecurve[showpoints=false,linewidth=2pt,linecolor=blue](3.4,8.2)(3.7,8.5)(3.9,8.9)
(3.2,9.3)(2.4,9)(2.1,8.2)(2.2,7.6)(2.6,7)(2.8,6.8)

\psecurve[showpoints=false,linewidth=2pt,linecolor=red](3.5,9.1)(3.8,9.2)(5,8.6)
(5.1,7.6)(4.6,6.9)(3.8,6.4)(3.4,6.2)
\psecurve[showpoints=false,linewidth=2pt,linecolor=red](3.8,6.4)(3.4,6.2)(2.9,6.1)
(2.6,6.2)(2.5,6.6)(2.8,7)(3.4,7.4)(3.7,7.5)
\psecurve[showpoints=false,linewidth=2pt,linecolor=red](3.4,7.4)(3.7,7.5)(3.9,7.9)
(3.3,8.8)(3.5,9.1)(3.8,9.2)

\rput(5.5,9.5){$6^2_2$}


\psecurve[showpoints=false,linewidth=2pt,linecolor=blue](9.9,8.7)(10.3,8.7)(11,8.8)
(11.6,8.6)(12,8.2)(11.8,7.1)(11.6,6.8)
\psecurve[showpoints=false,linewidth=2pt,linecolor=blue](11.8,7.1)(11.6,6.8)(11.2,6.3)
(10.8,6)(10.1,6.1)(9.8,6.4)
\psecurve[showpoints=false,linewidth=2pt,linecolor=blue](10.1,6.1)(9.8,6.4)(9.3,7)
(9.1,7.8)(9.3,8.5)(9.9,8.7)(10.3,8.7)

\psecurve[showpoints=false,linewidth=2pt,linecolor=red](11,6.3)(11.3,6.2)(11.9,6.6)
(11.7,7)(11.2,7.4)(10.6,7.7)(10.4,7.9)
\psecurve[showpoints=false,linewidth=2pt,linecolor=red](10.6,7.7)(10.4,7.9)(10.1,8.7)
(10.6,9.3)(11.1,8.9)(11.1,8.6)
\psecurve[showpoints=false,linewidth=2pt,linecolor=red](11.1,8.9)(11.1,8.6)(10.9,8.2)
(10,7.4)(9.5,7)(9.3,6.8)
\psecurve[showpoints=false,linewidth=2pt,linecolor=red](9.5,7)(9.3,6.8)(9.2,6.4)
(9.5,6)(10.4,6.6)(11,6.3)(11.3,6.2)

\rput(12.5,9.5){$7^2_1$}


\psecurve[showpoints=false,linewidth=2pt,linecolor=blue](2.6,2)(2.8,1.8)(3.4,1.4)
(3.8,1.2)(4.6,1.2)(4.7,1.8)(4.4,2)
\psecurve[showpoints=false,linewidth=2pt,linecolor=blue](4.7,1.8)(4.4,2)(4,2.2)
(3.6,2.4)(3.3,2.8)(3.4,3.2)(3.7,3.5)
\psecurve[showpoints=false,linewidth=2pt,linecolor=blue](3.4,3.2)(3.7,3.5)(3.9,3.9)
(3.2,4.3)(2.4,4)(2.1,3.2)(2.2,2.6)(2.6,2)(2.8,1.8)

\psecurve[showpoints=false,linewidth=2pt,linecolor=red](3.8,1.4)(3.4,1.2)(2.9,1.1)
(2.6,1.2)(2.5,1.6)(2.8,2)(3.4,2.4)(3.7,2.5)
\psecurve[showpoints=false,linewidth=2pt,linecolor=red](3.4,2.4)(3.7,2.5)(4.65,3.2)
(4.8,3.8)(4.6,4.2)(3.9,4.2)(3.5,4.1)
\psecurve[showpoints=false,linewidth=2pt,linecolor=red](3.9,4.2)(3.5,4.1)(3.3,4)
(3.2,3.7)(3.6,3.3)(4.4,3.1)(4.7,3)
\psecurve[showpoints=false,linewidth=2pt,linecolor=red](4.4,3.1)(4.7,3)(4.9,2.6)
(4.8,2.2)(4.4,1.8)(3.8,1.4)(3.4,1.2)

\rput(5.5,4.5){$7^2_2$}


\psecurve[showpoints=false,linewidth=2pt,linecolor=blue](10,3)(9.7,2.9)(9.4,2.9)
(9.2,3.4)(10,3.9)(11,3.9)(11.6,3.7)(11.9,3.4)(11.9,3.1)
\psecurve[showpoints=false,linewidth=2pt,linecolor=blue](11.9,3.4)(11.9,3.1)
(11.4,2.8)(10.6,3)(10,3)(9.7,2.9)

\psecurve[showpoints=false,linewidth=2pt,linecolor=red](9.3,3.3)(9.1,3)(9,2.4)
(9.4,1.5)(10,0.9)(10.4,0.8)(10.6,1.2)(10.5,1.5)(10.3,1.7)
\psecurve[showpoints=false,linewidth=2pt,linecolor=red](10.5,1.5)(10.3,1.7)
(10,2.3)(9.9,2.9)(9.7,3.3)(9.3,3.3)(9.1,3)

\psecurve[showpoints=false,linewidth=2pt,linecolor=green](10.3,1)(10.6,0.8)
(11.3,1)(12,2.2)(12,3.2)(11.6,3.3)(11.3,3)(11.2,2.7)
\psecurve[showpoints=false,linewidth=2pt,linecolor=green](11.3,3)(11.2,2.7)
(11,2.2)(10.6,1.8)(10.2,1.4)(10.3,1)(10.6,0.8)

\rput(12.5,4.5){$6^3_1$}


\endpspicture
\end{center}

\begin{center}
\rm Fig.~3
\end{center}
\end{exem}


\begin{thebibliography}{10}

\bibitem{ammann99a}
{\scshape B.~Ammann} -- {\og The {D}irac operator on collapsing
  {$S^1$}-bundles\fg}, In: Seminaire de theorie spectrale et geometrie. annee
  1997--1998, St. Martin D'Heres: Universite de Grenoble I, Institut Fourier,
  1998.

\bibitem{ammann98a}
\bysame , \emph{Spin-{S}trukturen und das {S}pektrum des {D}irac-{O}perators},
  Shaker Verlag, Aachen, 1998.

\bibitem{ammann00ppa}
\bysame , {\og Spectral estimates on 2-tori\fg}, Preprint, CUNY, Graduate
  Center, 2000.

\bibitem{ammann-baer98a}
{\scshape B.~Ammann {\normalfont \smfandname} C.~B{\"a}r} -- {\og The {D}irac
  operator on nilmanifolds and collapsing circle bundles\fg}, \emph{Ann. Glob.
  Anal. Geom.} \textbf{16} (1998), p.~221--253.

\bibitem{atiyah-patodi-singer75a}
{\scshape M.~F. Atiyah, V.~K. Patodi {\normalfont \smfandname} I.~M. Singer} --
  {\og Spectral asymmetry and {R}iemannian geometry {I}\fg}, \emph{Math. Proc.
  Camb. Phil. Soc.} \textbf{77} (1975), p.~43--69.

\bibitem{baer96a}
{\scshape C.~B{\"a}r} -- {\og The {D}irac operator on space forms of positive
  curvature\fg}, \emph{J. Math. Soc. Japan} \textbf{48} (1996), p.~69--83.

\bibitem{baer98ppb}
\bysame , {\og The {D}irac operator on hyperbolic manifolds of finite
  volume\fg}, {SFB256-Preprint} no. 566, Universit{\"a}t Bonn, 1998.

\bibitem{baer-schmutz92a}
{\scshape C.~B{\"a}r {\normalfont \smfandname} P.~Schmutz} -- {\og Harmonic
  spinors on {R}iemann surfaces\fg}, \emph{Ann. Glob. Anal. Geom.} \textbf{10}
  (1992), p.~263--273.

\bibitem{berline-getzler-vergne91a}
{\scshape N.~Berline, E.~Getzler {\normalfont \smfandname} M.~Vergne} --
  \emph{Heat kernels and {D}irac operators}, Springer-Verlag, Berlin
  Heidelberg, 1991.

\bibitem{bures93a}
{\scshape J.~Bures} -- {\og Spin structures and harmonic spinors on {R}iemann
  surfaces\fg}, In: Brackx, {F}. (ed.) et al., Clifford algebras and their
  applications in mathematical physics, Dordrecht: Kluwer Academic Publishers,
  1993.

\bibitem{bures98a}
\bysame , {\og Spin structures and harmonic spinors on nonhyperelliptc
  {R}iemann surfaces of small genera\fg}, In: Dietrich, {V}olker (ed.) et al.,
  Clifford algebras and their applications in mathematical physics, Dordrecht:
  Kluwer Academic Publishers, 1998.

\bibitem{chavel84a}
{\scshape I.~Chavel} -- \emph{Eigenvalues in {R}iemannian geometry}, Academic
  Press, Orlando etc., 1984.

\bibitem{dahl99ppa}
{\scshape M.~Dahl} -- {\og Dependence on the spin structure of the eta and
  {R}okhlin invariants\fg}, Preprint, Royal Insitute of Technology, Stockholm,
  1999.

\bibitem{friedrich84a}
{\scshape T.~Friedrich} -- {\og Zur {A}bh{\"a}ngigkeit des {D}irac-{O}perators
  von der {S}pin-{S}truktur\fg}, \emph{Coll. Math.} \textbf{48} (1984),
  p.~57--62.

\bibitem{gilkey84a}
{\scshape P.~Gilkey} -- \emph{Invariance theory, the heat equation and the
  {A}tiyah-{S}inger index theorem}, Publish or Perish, Wilmington, Delaware,
  1984.

\bibitem{gilkey89a}
\bysame , \emph{The geometry of sperical space form groups}, World Scientific,
  Singapore, 1989.

\bibitem{hitchin74a}
{\scshape N.~Hitchin} -- {\og Harmonic spinors\fg}, \emph{Adv. Math.}
  \textbf{14} (1974), p.~1--55.

\bibitem{lawson-michelsohn89a}
{\scshape H.~B. Lawson {\normalfont \smfandname} M.-L. Michelsohn} --
  \emph{Spin geometry}, Princeton University Press, Princeton, 1989.

\bibitem{lott00ppa}
{\scshape J.~Lott} -- {\og Collapsing and {D}irac-type operators\fg},
  \emph{ArXiv: math.DG/0005009} (2000).

\bibitem{pfaeffle99ppa}
{\scshape F.~Pf{\"a}ffle} -- {\og The {D}irac spectrum of {B}ieberbach
  manifolds\fg}, Preprint 1999, to app. in {J}. {G}eom. {P}hys.,
  Universit{\"a}t Hamburg, 1999.

\bibitem{sulanke79a}
{\scshape S.~Sulanke} -- {\og Berechnung des {S}pektrums des {Q}uadrates des
  {D}irac-{O}perators auf der {S}ph{\"a}re\fg}, Doktorarbeit, HU Berlin, 1979.

\bibitem{trautman95a}
{\scshape A.~Trautman} -- {\og The {D}irac operator on hypersurfaces\fg},
  \emph{Acta Phys. Polon. B} \textbf{26} (1995), p.~1283--1310.

\end{thebibliography}
\providecommand{\bysame}{\leavevmode ---\ }
\providecommand{\og}{``}
\providecommand{\fg}{''}
\providecommand{\smfandname}{et}
\providecommand{\smfedsname}{\'eds.}
\providecommand{\smfedname}{\'ed.}
\providecommand{\smfmastersthesisname}{M\'emoire}
\providecommand{\smfphdthesisname}{Th\`ese}

\end{document}